\documentclass{article}
\usepackage{amsfonts,latexsym}
\usepackage{amsmath, times}
\usepackage{amssymb}
\usepackage{amsthm}
\usepackage{color}
\usepackage[colorlinks=true, citecolor=blue, linkcolor=blue, urlcolor=black]{hyperref}

\newcommand{\R}{\mathbb R}

\newcommand{\del}{\partial}
\newcommand{\e}{\varepsilon}
\newtheorem{theorem}{Theorem}[section]
\newtheorem{lemma}{Lemma}

\newtheorem{corollary}{Corollary}[section]
\newtheorem{definition}{Definition}[section]
\newtheorem{remark}{Remark}[section]
\newlength{\defbaselineskip}
\newcommand{\setlinespacing}[2]%
{\setlength{\baselineskip}{#1 \defbaselineskip}}

\makeatletter

\@addtoreset{equation}{section} \makeatother 
\begin{document}
\begin{center}
 {\Large   {Ground state solutions for  weighted fourth-order Kirchhoff problem via Nehari method
}}
\end{center}
\vspace{0.2cm}
\begin{center}
Dridi Brahim   $^{(1)}$Rached Jaidane  $^{(2)}$  and   Rima Chetouane $^{(3)}$\\
  \
\noindent\footnotesize  $^{(1)}$ Umm Al-Qura University, Faculty of Applied Sciences, Department of mathematics, P.O. Box $14035$,
Holly Makkah $21955$, Saudi Arabia.\\
 Address e-mail: iodridi@uqu.edu.sa\\

 \
\noindent\footnotesize $^{(2)}$ Department of Mathematics, Faculty of Science of Tunis, University of Tunis El Manar, Tunisia.\\
 Address e-mail: rachedjaidane@gmail.com\\
 \noindent\footnotesize   $^{(3)}$ Department of Mathematics, Faculty of Exact Sciences, Frères Mentouri  Constantine 1 University, Algeria.
Address e-mail: rima.chetouane@umc.edu.dz
\end{center}

\vspace{0.5cm}
\noindent {\bf Abstract.}
In this article, we study the  following non local problem
$$g\big(\int_{B}w(x) |\Delta u|^{2}\big)\Delta(w(x)\Delta u) =|u|^{q-2}u +\ f(x,u) \quad\mbox{ in }\quad B, \quad u=\frac{\partial u}{\partial n}=0 \quad\mbox{ on } \quad\partial B,$$
 where $B$ is the unit ball in $\mathbb{R}^{4}$ and  $ w(x)$ is
a singular weight of logarithm type. The non-linearity is a combination of a reaction source
$f(x,u)$ which   is  critical in view of exponential inequality of Adams' type and  a polynomial function.
The Kirchhoff function $g$ is positive and continuous. By using the Nehari manifold method , the quantitative deformation lemma and degree theory results, we
establish the existence of a ground state solution. \\

\noindent {\footnotesize\emph{Keywords:} Weighted Sobolev space, Kirchhoff, biharmonic operator, Critical exponential growth, Adams' inequality.\\
\noindent {\bf $2010$ Mathematics Subject classification}: $35$J$20$, $49$J$45$, $35$K$57$, $35$J$60$.}

\section{Introduction and Main results}
In  this paper, we consider  the non local fourth order weighted elliptic equation:
\begin{equation}\label{eq:1.1}
   \displaystyle (P)~~~~\left\{
      \begin{array}{rclll}
g\big(\int_{B}w(x) |\Delta u|^{2}\big)\Delta(w(x)\Delta u)  &=& |u|^{q-2}u+ f(x,u)& \mbox{in} & B \\\\
        u&=& \displaystyle\frac{\partial u}{\partial n}=0 &\mbox{on }&  \partial B,
      \end{array}
    \right.
\end{equation}
where $B=B(0,1)$ is the unit open ball in $\R^{4}$, $q>4$, $f(x,t)$ is a radial function with respect to $x$, the weight  $w(x)$ is given by \begin{equation}\label{eq:1.2}w(x)=\bigg(\log \frac{e}{|x|}\bigg)^{\beta},~~\beta\in(0,1)\cdot\end{equation}
The Kirchhoff function $g$ is positive, continuous and verifies some mild conditions.
\\

The study of Kirchhoff problems was initiated in 1883, when Kirchhoff \cite{KIR} studied the following
equation

\begin{equation}\label{Kir}
	\rho\frac{\partial^2u}{\partial t^2}-\Big(\frac{P_0}{h} +
	\frac{E}{2L}\int_0^L|\frac{\partial u}{\partial
		x}|^2dx\Big)\frac{\partial^2u}{\partial x^2} = f(x,u),
\end{equation}
where $\rho, P_0, h, E, L$ represent physical quantities. This model extends the classical D’Alembert
wave equation by considering the effects of the changes in the length of the strings during
the vibrations. We call \eqref{Kir} a nonlocal problem since the equation contains an integral over $[0, L]$ which makes the study of it interesting.\
After  Lions in his pioneering work \cite{JL} presented an abstract functional analysis framework to \eqref{Kir}. We mention that non-local problems also arise in other areas, for instance, biological systems where the function $u$ describes a process that depends on the average of itself ( for example, population density), see for instance \cite{AJCTF,AJC}  and its references.\\

 Second order Kirchhoff's classical equation  has been extensively studied. We refer to the work of Chipot \cite{CL,CR}, Corr\^{e}a et al \cite{KIR} and their references. \\

 For $g(t)=1$, $w(x)=1$ and without the polynomial term,  we obtain the biharmonic equation $$\Delta^{2} u =f(x,u)~~\mbox{in}~~\mathbb{R}^{N}~~\mbox{or in bounded domain  }~~\Omega,$$ with Dirichlet or Navier boundary condition. The last equation  have been the subject of extensive mathematical studies in recent years, see \cite{AES, AL} and the references therein. \\Also, the biharmonic equation in dimension $N>4$ \begin{equation*}\Delta^{2} u=f(x,u)~~\mbox{in}~~\Omega\subset \mathbb{R}^{N},\end{equation*}where the nonlinearity $f$ has subcritical and critical polynomial growth of power less than $\frac{N+4}{N-4}$, have been extensively studied \cite{BGT,EFJ,GGS}. \\
 Also, for $g(t)=1$ where the weight $w(x)$ is given by (\ref{eq:1.2}) and without the polynomial term, Dridi and Jaidane \cite{DJ} studied the following weighted biharmonic problem  \begin{align}\nonumber \displaystyle \left\{
      \begin{array}{rclll}
    L:=\Delta(w(x) \Delta u)-\Delta u+V(x)u &=&  \displaystyle f(x,u)& \mbox{in} & B \\
        u=\frac{\partial u}{\partial n}&=&0 &\mbox{on }&  \partial B,
      \end{array}
    \right.
\end{align} $B$ is the unitary disk in $\mathbb{R}^{4}$, $f(x,t)$ is continuous and positive in
$B\times \mathbb{R}$ and behaves like $\exp\{\alpha t^{\frac{2}{1-\beta}}\}$ as $ |t|\rightarrow+ \infty$, for some $\alpha >0$ uniformly with respect to $x\in B$.   The potential $V :\overline{B}\rightarrow \mathbb{R}$ is a positive continuous function and bounded away from zero in $B$. The authors used Pass mountain method and minimax technics to establish the existence of radial solution.\\ An analogous problem but of the Kirchhoff type has been treated \cite{JR} that is \begin{align}\nonumber
  \displaystyle \left\{
      \begin{array}{rclll}
g\big(\int_{B}(w(x)|\Delta u|^{2}+|\nabla u|^{2} +V(x)u^{2})dx\big)\big[\Delta(w(x) \Delta u)-\Delta u+V(x)u) \big] &=& \ f(x,u)& \mbox{in} & B \\\\
 u=\displaystyle\frac{\partial u}{\partial n}&=&0 &\mbox{on }&  \partial B ,
      \end{array}
    \right.
\end{align}
In the latter problem, the non-linearity is positive  and behaves like $\exp\{\alpha t^{\frac{2}{1-\beta}}\}$ as $ |t|\rightarrow+ \infty$, for some $\alpha >0$ uniformly with respect to $x\in B,$ whereas in our case it is not. In our case the solution can change sign.\\
 Recently, Meng et al. \cite{MZZ}, studied the following fourth order equation of Kirchhoff type  namely:$$\Delta^{2}u-(a+b\int_{\mathbb{R}^{N}}|\nabla u|^{2}dx)\Delta u +V(x)u=f(x,u)~~\mbox{in}~~\mathbb{R}^{N},~~u\in H^{2}(\mathbb{R}^{N}),$$ with concave-convexe nonlinearities. The authors prove that there are at least two positive solutions. They used the Nehari manifold, Ekeland variational principle and the theory of Lagrange multipliers.\\
   we obtain the biharmonic equation

 In the literature, the notion of critical exponential growth is related to the Trudinger-Moser inegalities. For bounded domains $\Omega\subset \mathbb{R^{N}}$, and in the Sobolev space$ W^{1,N}_{0}(\Omega)$ , these inequalities \cite{JMo,NST}  are given by  $$\displaystyle\sup_{\int_{\Omega} |\nabla u|^{N}\leq1}\int_{\Omega}e^{\alpha|u|^{\frac{N}{N-1}}}~dx<+\infty~~\mbox{if and only if}~~\alpha\leq \alpha_{N},$$
where $\alpha_{N}=\omega_{N-1}^{\frac{1}{N-1}}$  with $\omega_{N-1}$ is the area of the unit sphere $S^{N-1}$ in $\mathbb{R}^{N}$.\\These kinds of results have permitted the study of second-order problems with exponential growth nonlinearities in non-weighted Sobolev spaces. For example, we cite the following problem of Kirchhoff type in dimension $N=2$ studied by Figueiredo and Severo \cite{FS}
\begin{equation*}
  \displaystyle \left\{
      \begin{array}{rclll}
-m\big(\int_{B}|\nabla u|^{2}dx\big)\Delta u  &=& \ f(x,u)& \mbox{in} & \Omega \\
       u &>&0 &\mbox{in }& \Omega\\
        u&=&0 &\mbox{on }&  \partial \Omega ,
      \end{array}
    \right.
\end{equation*}
where $\Omega$ is a smooth bounded domain in $\mathbb{R}^{2}$, the nonlinearity $f$ behaves like $\exp(\alpha t^{2})~~\mbox{as}~~t\rightarrow+\infty$, for some $\alpha>0$. $m:(0,+\infty)\rightarrow (0,+\infty)$ is a continuous function satisfying some conditions.  The authors proved that this problem has a positive ground state solution. The existence result was proved by combining minimax techniques and Trudinger-Moser inequality.\\

Recently, Sitong Chen, Xianhua Tang and Jiuyang Wei \cite{CTW}, studied the last problem. They have developed  some new approaches to estimate precisely the minimax level of
the energy functional and prove the existence of Nehari-type ground-state solutions and nontrivial solutions for the above
problem. \\
It should be noted that, recently, the following nonhomogeneous Kirchhoff-Schr\"{o}dinger equation
\begin{equation*}
  \displaystyle \left\{
      \begin{array}{rclll}
-M(\int_{\mathbb{R}^{2}}|\nabla u|^{2}+\xi(|x|)u^{2}dx )(-\Delta u +\xi(|x|)u)&=& Q(x)k(u)+\varepsilon h(x),\\
        u(x)\rightarrow0 ~~\mbox{as}~~|x|\rightarrow+\infty,
      \end{array}
    \right.
\end{equation*}
has been studied in \cite{AAU}, where
 $\varepsilon$ is a positive parameter, $M:\mathbb{R^{+}}\rightarrow\mathbb{R^{+}},$ $\xi ,Q:(0.+\infty)\rightarrow\mathbb{R}$, are continuous
functions that satisfy some mild conditions.  The nonlinearity $k:\mathbb{R}\rightarrow\mathbb{R}$ is continuous and  behaves like $\exp(\alpha t^{2})~~\mbox{as}~~t\rightarrow+\infty$, for some $\alpha>0$. The authors proved the existence of at least two weak solutions for this equation by combining the Mountain Pass Theorem
and Ekeland's Variational Principle. \\

These Trudinger-Moser inequalities mentioned above, have been extended to logarithmic weighted  Sobolev spaces  of the first order.  More precisely, Calanchi and Ruff proved the following result  :\begin{theorem}\label{th1}\cite{CR2} \label{th1.1}\begin{itemize}\item[$(i)$] ~~Let $\beta\in[0,1)$ and let $\tau$ the weight given by $ \tau(x)=\big(\log \frac{1}{|x|}\big)^{\beta}$, then
  \begin{equation*}
 \int_{B} e^{|u|^{\xi}} ~dx <+\infty, ~~\forall~~u\in W_{0,rad}^{1,N}(B,\tau),~~
  \mbox{if and only if}~~\xi\leq \xi_{N,\beta}=\frac{N}{(N-1)(1-\beta)}=\frac{N'}{1-\beta}
 \end{equation*}
and
 \begin{equation*}
 \sup_{\substack{u\in W_{0,rad}^{1,N}(B,\tau) \\ \int_{B}|\nabla u|^{N}w(x)~dx\leq 1}}
 \int_{B}~e^{\alpha|u|^{\xi_{N,\beta} }}~dx < +\infty~~~~\Leftrightarrow~~~~ \alpha\leq \alpha_{N,\beta}=N[\omega^{\frac{1}{N-1}}_{N-1}(1-\beta)]^{\frac{1}{1-\beta}}
 \end{equation*}
where $\omega_{N-1}$ is the area of the unit sphere $S^{N-1}$ in $\R^{N}$ and $N'$ is the H$\ddot{o}$lder conjugate of $N$.
\item [$(ii)$] Let $\sigma$ the weight given by $\sigma(x)=\big(\log \frac{e}{|x|}\big)^{N-1}$, then
  \begin{equation*}\label{eq:71.5}
 \int_{B}exp\{e^{|u|^{\frac{N}{N-1}}}\}~dx <+\infty, ~~~~\forall~~u\in W_{0,rad}^{1,N}(B,\sigma)
 \end{equation*} and
 \begin{equation*}\label{eq:71.6}
\sup_{\substack{u\in W_{0,rad}^{1,N}(B,\sigma) \\  \|u\|_{\sigma}\leq 1}}
 \int_{B}exp\{\beta e^{\omega_{N-1}^{\frac{1}{N-1}}|u|^{\frac{N}{N-1}}}\}~dx < +\infty~~~~\Leftrightarrow~~~~ \beta\leq N,
 \end{equation*}
\end{itemize}
\end{theorem}
where $B$ is the unit ball of $\mathbb{R^{N}}$, $N\geq2$ and \newline$W_{0,rad}^{1,N}(B,\sigma)=\mbox{closure}\{u \in
C_{0,rad}^{\infty}(B)~~|~~\int_{B}\sigma(x)|\nabla u|^{N}~dx <\infty\}$, is equipped with the norm $$\displaystyle\|u\|_{\sigma}=\big(\int_{B}\sigma(x)|\nabla u|^{N}~dx\big)^{\frac{1}{N}}\cdot$$\\This has allowed the study of problems involving logarithmic weight operators and exponential growth nonlinearities. We cite
that  recently, in the case, $V=0$ or $V\neq 0$, Baraket et al. \cite{BJ}, S. Deng, T. Hu and C-L. Tang and Calanchi et al. \cite{CRS,DHT}, have proved the existence of a nontrivial solution for
the following boundary value problem
$$
 \displaystyle \left\{
      \begin{array}{rclll}
    -\textmd{div} (\sigma(x)|\nabla u(x)|^{N-2}\nabla u(x) ) +V(x)|u|^{N-2}u &=& \ f(x,u)& \mbox{in} & B \\
        u&=&0 &\mbox{on }&  \partial B,
      \end{array}
         \right.
      $$
where $B$ is the unit ball in $\R^N,\; N\geq2,$ the radial positive weight $\sigma(x)=(\log\frac{e}{|x|})^{N-1}$, the function $f(x, u)$ is continuous in $B\times\R$ and behaves like $\exp\big(e^{\alpha{t^{\frac{N}{N-1}}}}\big)~~\mbox{as}~~t\rightarrow+\infty$, for some $\alpha>0$. The authors proved that there is a non-trivial solution to this problem using minimax techniques combined with weighted Trudinger-Moser inequality. \\

Also, we point out that recently, in the case where $g$ is not constant, Abid et al. and Jaidane \cite{ABJ,J} have proved the existence of a nontrivial solution for
the following logarithmic weighted Kirchhoff problem
$$
 \displaystyle \left\{
      \begin{array}{rclll}
    -g\Big( \displaystyle\int_{B}\tau(x)|\nabla u|^{N}+\overline{V}(x)|u|^{N}dx\Big)\textmd{div} (\tau(x)|\nabla u|^{N-2} \nabla u+\overline{V}(x)|u|^{N-2}u) &=& \ f(x,u)& \mbox{in} & B \\
        u&=&0 &\mbox{on }&  \partial B,
      \end{array}
         \right.
      $$
where $B$ is the unit ball in $\mathbb{R}^{N},\; N\geq2,$ the weight $\tau(x)=\Big(\log \frac{e}{|x|}\Big)^{\beta(N-1)},~~\mbox{with}~~\beta=1~~\mbox{or}~~\beta\in[0,1),$
 the reaction term  $f(x, u)$ is continuous in $B\times \mathbb{R}$ and behaves like $\exp \big(e^{\alpha t^{\frac{N}{(N-1)}}}\big)~~\mbox{or}~~e^{\alpha t^{\frac{N}{(N-1)(1-\beta)}}},~~\mbox{as}~~t\rightarrow+\infty$, for some $\alpha>0$ and the potential $\overline{V}$ is a positive and continuous function on $\overline{B}$ or equal to zero. The authors proved that there is a non-trivial solution to this problem using Nehari method and weighted  Trudinger-Moser inequality \cite{CR2}.\\

We recall that the Trudinger-Moser inequalities have been extended to spaces of higher orders to finally obtain the so-called Adams' inequalities \cite{Ada,BRFS}. In fact, for bounded domains $\Omega\subset \mathbb{R}^{4}$, in \cite{Ada,BRFS} and in the space $ W_{0}^{2,2}(\Omega)$ the authors  obtained, \begin{equation*}
 \sup_{\substack{u\in S }}
 \int_{\Omega}~(e^{\alpha u^{2})}-1)dx < +\infty~~~~\Leftrightarrow~~~~ \alpha\leq 32 \pi^{2}
 \end{equation*}
 where $S=\{u\in W^{2,2}_{0}(\Omega)~~|~~\displaystyle\big(\int_{\Omega}|\Delta u|^{2}dx \big)^{\frac{1}{2}}\leq 1\}.$ \\This last result opened the way to study fourth-order problems with subcritical or critical nonlinearity. We cite the work of Sani \cite{FS} \begin{equation*}\Delta^{2} u+V(x) u=f(x,u)~~\mbox{in}~~H^{2}(\mathbb{R}^{4})\cdot\end{equation*}

Recently, a result related to Adams inequalities has been extended to Sobolev spaces with logarithmic weights. Lately, Wang and Zhu \cite{WZ} proved the following result:
 \begin{theorem}\cite{WZ} \label{th1.1} ~~Let $\beta\in(0,1)$ and let $w$ given by (\ref{eq:1.2}), then
\begin{equation}\label{eq:1.3}
 \sup_{\substack{u\in W_{0,rad}^{2,2}(B,w) \\  \int_{B}w(x)|\Delta u|^{2}~dx \leq 1}}
 \int_{B}~e^{\displaystyle\alpha|u|^{\frac{2}{1-\beta} }}~dx < +\infty~~~~\Leftrightarrow~~~~ \alpha\leq \alpha_{\beta}=4[8\pi^{2}(1-\beta)]^{\frac{1}{1-\beta}},
 \end{equation}

 \end{theorem}
 where $W_{0,rad}^{2,2}(B,w)$ denotes the weighted Sobolev space of radial functions given by $$W_{0,rad}^{2,2}(B,w)=\mbox{closure}\{u\in
C_{0,rad}^{\infty}(B)~~|~~\displaystyle\int_{B}w(x)|\Delta u|^{2}~dx <\infty\},$$ endowed with the norm  $\displaystyle\|u\|= \big(\int_{B}w(x)|\Delta u|^{2}~dx\big)^{\frac{1}{2}}$.\\

Before stating our results, let's start by defining our functional space. Now, let $\Omega \subset \R^{4}$, be a bounded domain and $w\in L^{1}(\Omega)$ be a nonnegative function. We introduce the Sobolev space$$ W_{0}^{2,2}(\Omega,w)=\mbox{closure}\{u\in
C_{0}^{\infty}(\Omega)~~|~~\displaystyle\int_{\Omega}w(x)|\Delta u|^{2}dx <\infty\}.$$
We will focus on radial functions and consider the naturel subspace
$$\mathbf{W}:= W_{0,rad}^{2,2}(B,w)=\mbox{closure}\{u\in
C_{0,rad}^{\infty}(B)~~|~~\displaystyle\int_{\Omega}w(x)|\Delta u|^{2}dx <\infty\},$$ with the specific weight $w$ given by (\ref{eq:1.2}).
The subspace $\mathbf{W}$ is endowed with the norm
\begin{equation*}
   \| u\|=\left(\int_{B}w(x)|\Delta u|^{{2}}dx\right )^{\frac{1}{2}}.
\end{equation*}
 We note that this norm is  induced from the scalar product,
$$\langle u,v\rangle=\int_{B}\big( w(x)\Delta u\Delta v \big)~dx.$$

Let $\gamma:=\displaystyle\frac{2}{1-\beta}$. According to  inequality (\ref{eq:1.3}), we will say that $f$ has critical growth at infinity if there exists some $\alpha_{0}>0$,
\begin{equation}\label{eq:1.6}
\lim_{s\rightarrow +\infty}\frac{|f(x,s)|}{e^{\alpha s^{\gamma}}}=0,~~~\forall~\alpha ~~\mbox{such that}~~ \alpha_{0}<\alpha~~
\mbox{and}~~~~\lim_{s\rightarrow +\infty}\frac{|f(x,s)|}{e^{\alpha s^{\gamma}}}=+\infty,~~\forall~\alpha<\alpha_{0}.\end{equation}

Now we define the Kirchhoff function $g$  and give the conditions on it . The function $g$ is continuously differentiable in $\mathbb{R^{+}}$ and verifies :
\begin{description}
	\item[$(G_{1})$]  $g$ is increasing with $g(0)=g_{0}>0$;
	\item[$(G_{2})$] $\displaystyle\frac{g(t)}{t} ~~\mbox{is  nonincreasing for} ~~ t>0.$
\end{description}The assmption $(G_{2})$ implies that $\displaystyle \frac{g(t)}{t} \leq g(1)$ for all $t\geq 1.$ \\
From $(G_{1})$ and $(G_{2})$, we can get
	\begin{equation}\label{eq:1.5} G(t+s)\geq G(t)+G(s)~~\forall ~~s,t\geq0~~ \mbox{where}~~G(t)=\int^{t}_{0}g(s)ds\end{equation}
and
	
\begin{equation}\label{eq:1.6}g(t)\leq g(1) + g(1) t,\quad \forall~~ t\geq 0\cdot\end{equation}
	As a consequence, we get \begin{equation}\label{eq:1.7}G(t)\leq g(1)t + \frac{g(1)}{2} t^{2},\quad \forall~~ t\geq 0\cdot\end{equation}
Moreover, we have that

\begin{equation}\label{eq:1.8}\frac{1}{2}G(t)-\frac{1}{4}g(t)t~~\mbox{is nondecreasing and positive for}~~t > 0.\end{equation}
Consequently, one has for all $q\geq4$, \begin{equation}\label{eq:1.7'}t\mapsto \frac{1}{2}G(t)-\frac{1}{q}g(t)t,~~\mbox{is nondecreasing and positive for}~~t > 0.\end{equation}
A typical example of a function $g$ fulfilling the conditions $(G_{1})$ and $(G_{2})$  is given by
$$g(t)=g_{0}+at,~~g_{0},a>0.$$
Another example is given by $g(t)=1+\ln (1+t)$.\\

Furthermore, we suppose that $f(x,t)$ has critical growth and satisfies the following hypothesis:
\begin{enumerate}
\item[$(A_{1})$] $f: B \times \mathbb{R}\rightarrow\mathbb{R}$ is continuous  and radial in $x$.
\item[$(A_{2})$] There exist $\theta > q\geq4$  such that  we have
 $$0 < \theta F(x,t)\leq tf(x,t), \forall (x, t)\in~~B\times\mathbb{R} \setminus\{0\} $$
  where
$$F(x,t)=\displaystyle\int_{0}^{t}f(x,s)ds.$$
\item [$(A_{3})$]  For each $x\in B$,~$\displaystyle t \mapsto \frac{f(x,t)}{|t|^{q-1}}~~\mbox{is increasing for}~~ t\in~~\mathbb{R} \setminus\{0\}$.
 \item [$(A_{4})$]   $\displaystyle\lim_{t\rightarrow 0}\frac{|f(x,t)|}{|t|}=0.$
 \item[$(A_{5})$] There exist $p$ such that $p>q>4 $ and $C_p > 1$ such that
$$ sgn(t) f(x,t) \geq C_p \vert t\vert^{p-1}, \quad \mbox{for all} \; (x,t)\in B \times \R,$$
where $sgn(t) = 1$ if $t > 0,$ $sgn(t) = 0$ if $t = 0,$ and $sgn(t) = -1$ if $t < 0.$
\end{enumerate}

\begin{remark}\label{rem1}The conditions $(A_{2})$ and $(A_{3})$ implies that $\displaystyle t\mapsto \frac{f(x,t)}{t^{3}}~~\mbox{is increasing for}~~ t>0$.
\end{remark}
We give an example of such nonlinearity. The nonlinearity $f(x,t)=C_{p}|t|^{p-2}t+|t|^{p-2}t\exp(\alpha_{0} |t|^{\gamma})$ satisfies the assumptions $(A_{1})$, $(A_{2})$, $(A_{3})$ ,$(A_{4})$ and $(A_{5})$.\\

We will consider the following definition of  solutions.
\begin{definition}\label {def1.1}
We say that a function $u\in \mathbf{W}$ is a weak solution to problem $(P)$ if
\begin{equation*}\label {eq:1.9}
g(\|u\|^{2})\int_{B}w(x)~\Delta u ~\Delta \varphi~  dx =\int_{B}|u|^{q-2}u\varphi~~dx+
\int_{B}f(x,u)~~ \varphi ~dx,~~\forall~\varphi \in \mathbf{W}.
\end{equation*}
\end{definition}
Let $\mathcal{J} :\mathbf{W} \rightarrow \R$ be the energy functional given by
 \begin{equation}\label{eq:1.5}
\mathcal{J}(u)=\frac{1}{2}G(\|u\|^{2})-\frac{1}{q}\int_{B}|u|^{q}~dx -\int_{B}F(x,u)~dx,
\end{equation}
where
$$ F(x,t)=\displaystyle\int_{0}^{t}f(x,s)ds.$$
The energy $\mathcal{J}~~\mbox{belongs to}~~  C^{1}(\mathbf{W},\mathbb{R})$ (see Section 3)  and $$\langle\mathcal{J}'(u),\varphi\rangle:=\mathcal{J'}(u)\varphi=g(\|u\|^{2})\langle u,\varphi\rangle-\int_{B}|u|^{q-2}u\varphi~~dx-\int_{B}f(x,u)~ \varphi~dx~,~~\forall~~\varphi \in\mathbf{W}\cdot$$
Our strategy is to identify solutions that minimize the associated energy $\mathcal{J}$ from all solutions to problem $(P)$. For this goal, we define the Nehari set as follows

$$\displaystyle\ \mathcal{N}:=\{u\in \mathbf{W}:\langle \mathcal{J}'(u),u\rangle=0, u\neq 0\}$$
and we are looking for a minimization of the energy function $\mathcal{J}$  through the following minimization problem:
\begin{equation*}\label {eq:1.12}
\displaystyle m=\inf_{u\in\mathcal{N}}\mathcal{J}(u)\cdot
\end{equation*}

To our best knowledge, there are no results for  solutions to the non local weighted biharmonic  equation with critical exponential nonlinearity combined with a polynomial term on the weighted Sobolev space $\mathbf{W}$.\\

Now, we give our main result as follows:
\begin{theorem}\label{th1.3}~~ Let $f(x,t)$ be a function that has a critical growth at
$+\infty$. Suppose that   $(A_{1})$, $(A_{2})$, $(A_{3})$, $(A_{4})$ , $(A_{5})$,  $(G_{1})$ and  $(G_{2})$  are satisfied.
 Then  problem $(P)$ has a radial solution with minimal energy provided
\begin{equation}\label{eq:1.9} C_{p} >\max\biggm\{1, 2\tau^{\frac{p}{2}}\bigg(\frac{4q^{2}(p-2)m_{p}}{g_{0}(q-4)(p-q)}\big(\frac{2(\alpha_{0}+\delta)}{\alpha_{\beta}}\big)^{1-\beta}\bigg)^{\frac{p-2}{2}}\biggm\}
\end{equation}
 where $\delta >0$, $\displaystyle\tau=\frac{g(1)}{2g_{0}}+\frac{g(1)}{g^{2}_{0}}\frac{p}{p-4}m_{p}$~~ \mbox{with}~~ $m_p = \displaystyle\inf_{u\in\mathcal{N}_p} J_p(u)>0$,
 $$J_p(u) :=\frac{1}{2}G(\|u\|^{2}) -\frac{1}{p}\int_{B} \vert u\vert^p dx  $$ and
 $$\mathcal{N}_p:= \{ u \in \mathbf{W}, u \neq 0\; \mbox{and}\; \langle J_p'(u),u\rangle =0\}.$$
 \end{theorem}

In general the study of  fourth order partial differential equations is considered an interesting topic. The interest in studying such equations was stimulated by their applications in micro-electro-mechanical systems, phase field models of multi-phase systems, thin film theory, surface diffusion on solids, interface dynamics, flow in Hele-Shaw cells, see \cite{D, FW, M}.\\

This present work is organized as follows: in Section $2$, we present some necessary
preliminary knowledge about functional space and some preliminaries results. In Section $3$, we introduce the variational
framework and some technical key lemmas . In Section $4$, we study an auxiliary problem which will be of great use to prove our main result. Section $5$ is devoted to the proof of Theorem \ref{th1.3}.\\
Finally,  we note that a constant $C$ may change from one line to another and sometimes we index the constants in order to show how they change. Also, we shall use the notation $\vert u \vert_p$ for the norm in the Lesbegue space $L^p(B)$.

\section{Weighted Sobolev Space setting and preliminaries results }
Let $\Omega \subset \R^{N}$, $N\geq2$,  be a bounded domain in $\R^{N}$ and let $w\in L^{1}(\Omega)$ be a nonnegative function. To handle with weighted operator, we need to introduce some functional spaces $L^{p}(\Omega,w)$, $W^{m,p}(\Omega,w)$, $W_{0}^{m,p}(\Omega,w)$. To this end, let $S(\Omega)$ be the set of all measurable real-valued functions defined on $\Omega$ and two measurable functions are considered as the same element if they are equal almost everywhere.\\

According to  Drabek et al. and Kufner in \cite{DKN,Kuf}, the weighted Lebesgue space $L^{p}(\Omega,w)$ is set as follows:
$$L^{p}(\Omega,w)=\{u:\Omega\rightarrow \R ~\mbox{measurable};~~\int_{\Omega}w(x)|u|^{p}~dx<\infty\}$$
for any real number $1\leq p<\infty$.\\
This is a normed vector space equipped with the norm
$$\|u\|_{p,w}=\Big(\int_{\Omega}w(x)|u|^{p}~dx\Big)^{\frac{1}{p}}.$$
For $w(x)=1$, one finds the standard  Lebesgue space $L^{p}(\Omega)$ endowed with the norm $|u|_{p}=\Big(\int_{\Omega}|u|^{p}~dx\Big)^{\frac{1}{p}}.$\\

For $m\geq 2$, let $w$ be a given family of weight functions $w_{\tau}, ~~|\tau|\leq m,$ $$w=\{w_{\tau}(x)~~x\in\Omega,~~|\tau|\leq m\}.$$
In \cite{DKN}, the  corresponding weighted Sobolev space was  defined as
$$ W^{m,p}(\Omega,w)=\{ u \in L^{p}(\Omega)~ ~\mbox{such that}~~ D^{\tau} u \in L^{p}(\Omega)~\mbox{for all}~~|\tau|\leq m-1,~~D^{\tau} u \in L^{p}(\Omega,w)~\mbox{for all}~~|\tau|= m\} $$
endowed with the following norm:

\begin{equation*}\label{eq:2.2}
\|u\|_{W^{m,p}(\Omega,w)}=\bigg(\sum_{ |\tau|\leq m-1}\int_{\Omega}|D^{\tau}u|^{p}dx+\displaystyle \sum_{ |\tau|= m}\int_{\Omega}w(x)|D^{\tau}u|^{p} dx\bigg)^{\frac{1}{p}}.
\end{equation*}

If we also assume that $w(x)\in L^{1}_{loc}(\Omega)$, then $C^{\infty}_{0}(\Omega)$ is a subset of $W^{m,p}(\Omega,w)$. Hence, we can introduce the space $W^{m,p}_{0}(\Omega,w)$,
as the closure of $C^{\infty}_{0}(\Omega)$ in $W^{m,p}(\Omega,w).$ Moreover, the following embedding is compact  $$W^{m,p}(\Omega,w)\hookrightarrow\hookrightarrow W^{m-1,p}(\Omega)\cdot$$
Also, $(L^{p}(\Omega,w),\|\cdot\|_{p,w})$ and $(W^{m,p}(\Omega,w),\|\cdot\|_{W^{m,p}(\Omega,w)})$ are separable, reflexive Banach spaces provided that $w(x)^{\frac{-1}{p-1}} \in L^{1}_{loc}(\Omega)$.\\
 Then the space $\mathbf{W}$ is a Banach and reflexive space with the norm
$$ \| u\|=\left(\int_{B}w(x)|\Delta u|^{{2}}dx \right )^{\frac{1}{2}}$$
which is equivalent to the following norm (see lemma \ref{lemr} bellow) $$\|u\|_{W_{0,rad}^{2,2}(B,w)}=\displaystyle\big(\int_{B}u^{2}dx +\int_{B} |\nabla u|^{2}~dx+\int_{B}w(x)|\Delta u|^{2}dx \big)^{\frac{1}{2}}\cdot$$
We also have the continuous embedding   $$\mathbf{W}\hookrightarrow L^{q}(B)~~\mbox{for all}~~q\geq 1.$$
 Moreover, $\mathbf{W}$ is compactly
embedded in $L^{q}(B)$  for all $q \geq 2$ . In fact, we have



\begin{lemma}\label{lemr}

\item[(i)] Let $u$ be a radially symmetric
 function in $C_{0}^{2}(B)$. Then, we have\begin{itemize}\item[$(i)$]\cite{WZ}
 $$|u(x)|\leq \displaystyle\frac{1}{2\sqrt{2}\pi}\frac{||\log(\frac{e}{|x|})|^{1-\beta}-1|^{\frac{1}{2}}}{\sqrt{1-\beta}}\displaystyle \int_B w(x)|\Delta u|^2dx \leq \displaystyle\frac{1}{2\sqrt{2}\pi}\frac{||\log(\frac{e}{|x|})|^{1-\beta}-1|^{\frac{1}{2}}}{\sqrt{1-\beta}}\|u\|^{2}\cdot$$

\item[(ii)] The norms $\|.\| $ and $\|u\|_{ W_{0,rad}^{2,2}(B,w)}=\displaystyle\big(\int_{B}u^{2}dx +\int_{B} |\nabla u|^{2}~dx+\int_{B}w(x)|\Delta u|^{2}dx \big)^{\frac{1}{2}}$ are equivalents.
\item[(iii)]  The following embedding
is continuous $$\mathbf{W}\hookrightarrow L^{q}(B)~~\mbox{for all}~~q\geq 1.$$
\item[(iv)]$\mathbf{W}$ is compactly
embedded in $L^{q}(B)$  for all $q \geq2$ .
\end{itemize}
\end{lemma}
\textit{Proof }

$(i)$ see \cite{WZ}\\

 $(ii)$~~ By
 Poincar\'{e} inequality,  for all $u \in W_{0,rad}^{1,2}(B)$
 $$\int_B  \displaystyle  u^{2}~dx
 \leq C \int_B  \displaystyle |\nabla u|^{2}~dx.
 $$ Using the Green formula, we get
$$
\int_B  \displaystyle |\nabla u|^{2}~dx= \int_B \nabla u. \nabla u ~dx = - \displaystyle \int_B u \Delta
 u ~dx + \underbrace{\displaystyle\int_{\del B} u \frac{\partial u}{\partial n}}_{= 0}dn ~
 \leq  \displaystyle \Big| \int_B u \Delta u ~dx\Big|\cdot
 $$
By  Young inequality, we get for all $\varepsilon>0$
 $$
 \displaystyle \Big| \int_B u \Delta u \Big|~dx \leq
  \displaystyle \frac{1}{2 \varepsilon} \displaystyle \int_B |\Delta u|^2 ~dx + \displaystyle \frac{\varepsilon}{2} \displaystyle
 \int_B u^2~dx\leq  \displaystyle \frac{1}{2 \varepsilon} \displaystyle \int_B w(x)|\Delta u|^2 ~dx + \displaystyle \frac{\varepsilon}{2} \displaystyle
 \int_B u^2 ~dx.$$
 Hence
 $$ (1 - \displaystyle \frac{\varepsilon}{2} C) \displaystyle \int_B |\nabla u |^{2} ~dx~ \leq ~
 \displaystyle \frac{1}{2\varepsilon} \displaystyle\int_ B  w(x)|\Delta u |^{2}~dx,$$
  then, $$\displaystyle\int_{B}u^{2}~dx +\int_{B} |\nabla u|^{2}~dx+\int_{B} w(x)|\Delta u|^{2}~dx \leq C \int_{B} w(x)|\Delta u|^{2}~dx\leq C\|u\|^{2}\cdot$$
 Then $(ii)$ follows.\\
 $(iii)$ and $(iv)$. Since $w(x)\geq1$, then following embedding are continuous and compact   $$\mathbf{W}\hookrightarrow W^{2,2}_{0,rad}(B,w)\hookrightarrow W^{2,2}_{0,rad}(B)\hookrightarrow L^{q}(B)~~\forall q\geq 2$$ and from $(i)$, we have that $\mathbf{W} \hookrightarrow L^{q}(B)$ is continuous for all $q\geq1$.
This concludes the lemma.\hfill $\Box$\\

\section{ The variational setting and some technical lemmas }
Note that, by the hypothesis ($A_{4}$), for any $\varepsilon>0$, there exists $\delta_{0}>0$
   such that \begin{equation}\label{e1}|f(x,t)|\leq \varepsilon |t|,~~\forall ~~0<|t|\leq \delta_{0} ,~~\mbox{uniformly in}~~ x\in B.
   \end{equation}
   Moreover, since $f$ is critical at infinity, for every $\varepsilon>0$, there exists $C_{\varepsilon}>0 $ such that
   \begin{equation}\label{e2}~~\forall t\geq C_{\varepsilon}~~|f(x,t)|\leq \varepsilon \exp( ~c|t|^{\gamma})~~\mbox{with}~~c>\alpha_{0}~~\mbox{uniformly in}~~ x\in B.
   \end{equation}In particular, we obtain for $q>2$,\begin{equation}\label{e3} ~|f(x,t)t|\leq \frac{\varepsilon}{C^{q-1}_{\varepsilon}}|t|^{q} \exp(c ~|t|^{\gamma})~~\mbox{with}~~c>\alpha_{0}~~\mbox{uniformly in}~~ x\in B.
   \end{equation}
 Hence, using   (\ref{e1}), (\ref{e2}), (\ref{e3})  and the continuity of $f$,  for every $\varepsilon>0$, for every $q>2$, there exist positive constants $C$ and $c$ such that \begin{align}\label{imp} |f(x,t)|\leq \varepsilon |t| +C |t|^{q-1}e^{c ~|t|^{\gamma}}, ~~~~~~\forall\  (x,t)\in B\times \mathbb{R}.\end{align}
 It follows from (\ref{imp}) and $(A_{2})$, that for all $\varepsilon>0$, there exists $C>0$ such that
\begin{equation}\label {eq:1.10}
F(x,t)\leq \frac{1}{4}\varepsilon|t|^{2}+C |t|^{q}e^{a~|t|^{\gamma}},~~~~~~\mbox{for all}~~t\in \R.
\end{equation}

 So, by (\ref{eq:1.5}) and  (\ref{eq:1.10}) the functional $\mathcal{J}$ given by (\ref{eq:1.5}), is well defined. Moreover, by standard arguments,  $\mathcal{J}\in  C^{1}(\mathbf{W},\mathbb{R})$.
It is standard to check that critical points of $\mathcal{J}$ are precisely weak solutions of $(P)$. Moreover, we have $$\langle\mathcal{J}'(u),\varphi\rangle=\mathcal{J'}(u)\varphi=g(\|u\|^{2})\int_{B}\big(w(x)~ \Delta u~\Delta \varphi\big)~dx-\int_{B}|u|^{q-2}u\varphi~~dx-\int_{B}f(x,u)~ \varphi~dx~,~~\forall~~\varphi \in\mathbf{W}\cdot$$

In the next results, we show that $\mathcal{N}$ is not empty and that $\mathcal{J}$, restricted to $\mathcal{N}$, is bounded from below.\\
In the following we assume, unless otherwise stated, that the function $f$ satisfies the conditions $(A_{1})$ to $(A_{4})$ and the function $g$ satisfies $(G_{1})$ and $(G_{2})$ . Let $u\in\mathbf{W}$ with $u\not\equiv 0 $ a.e. in the ball $ B$, and we deﬁne the function
$\displaystyle\Upsilon_{u} : [0, \infty)  \rightarrow\R$
 as
\begin{equation}\label{eq:3.1}
 \displaystyle\Upsilon_{u}(t) = \mathcal{J}(tu ).
\end{equation}
It's clear that $\displaystyle\Upsilon'_{u}(t)=0$ is equivalent to $tu\in \mathcal{N}$.
\begin{lemma}\label{lem1}\begin{itemize}
\item [(i)]For  each  $u \in\mathbf{W}$ with  $u\neq0$  ,  there  exists  an  unique $t_{u}>0$,
 such that
$\displaystyle t_{u}u \in\mathcal{N}.$ In particular,  the set $\mathcal{N}$ is nonempty and $\mathcal{J}(u)>0$, for every $u\in \mathcal{N}$.
\item[(ii)]  For all $t \geq 0 $ with $t \neq t_{u},$ we have
$$\displaystyle \mathcal{J}(tu ) <  \mathcal{J}(t_{u}u)\cdot$$\end{itemize}

\end{lemma}
Proof.$(i)$ Note that since $q\geq 4$, we have
$$\displaystyle\lim_{|t|\rightarrow 0} \frac{|t|^{q-1}}{|t|}=0,$$
$$\displaystyle\lim_{|t|\rightarrow \infty} \frac{|t|^{q-1}}{|t|^{r-1}}=0, \mbox{  for  all } r\in (q, \infty),$$
  Then for any  $\epsilon > 0$, there exists a positive constant $C_{1} = C_{1}(\varepsilon )$ such that
\begin{equation}\label{eq3}
 |t|^{q-1} \leq \epsilon |t|+ C_{1}|t|^{r-1}\mbox{  for  all }     t \in \R.
\end{equation}\\

From (\ref{imp}) and (\ref{eq:1.10}), for all $\epsilon> 0$, there exist  positive constants $C'_{1} = C_{1}(\epsilon )$  and $C_{2} = C_{2}(\epsilon )$ such that
\begin{equation}\label{eq:3.2}
f(x,t)t\leq \epsilon |t|^{2} +C'_{1}|t|^{r }\exp(\alpha|t|^{\gamma}  )\mbox{ for all }	\alpha > \alpha_{0}, r >q.
\end{equation}
and \begin{equation}\label{eq:3.3}
F(x,t)\leq  \frac{1}{4}\epsilon |t|^{2} +C_{2}|t|^{r }\exp(\alpha|t|^{\gamma}  )\mbox{ for all }	\alpha > \alpha_{0}, r >q .
\end{equation}

 Now, given $u\in\mathbf{W}$ fixed with  $u \neq0$ . From (\ref{eq:3.3}), (\ref{eq3}) and (\ref{eq:1.8}), for all $\varepsilon>0$, we have
\begin{align}\label{eq:3.4}
 \displaystyle\Upsilon_{u}(t) = \mathcal{J}(tu )
 & =\frac{1}{2}G(t^{2}\|u\|^{2})-\frac{1}{q}\int_{B}|tu|^{q}~dx-\int_{B}F(x,tu)~dx\nonumber   \\
 & \geq \frac{1}{4}g(t^{2}\|u\|^{2})t^{2}\|u\|^{2}-\frac{\epsilon}{q} t^{2} \int_{B}|u|^{2}~dx -C_{1}\frac{|t|^{r}}{q}\int_{B}|u|^{r}-\int_{B}F(x,tu)~dx\nonumber \\
 & \geq \frac{g_{0}}{4}t^{2}\|u\|^{2} \nonumber
 -\frac{1}{4}\epsilon t^{2} \int_{B}|u|^{2}dx\nonumber-\frac{\epsilon}{q} t^{2} \int_{B}|u|^{2}~dx\\ &-C_{1}\frac{|t|^{r}}{q}\int_{B}|u|^{r} - C'_{1}\int_{B}|tu|^{r}\exp (\alpha t|u|^{\gamma})dx
 \end{align}
   Using the
 H\"{o}lder inequality, with $a, a' > 1$ such that $\displaystyle\frac{1}{a} + \frac{1 }{a'}  = 1$, and Sobolev embedding Lemma \ref{lemr}, we get
   \begin{align}
 \displaystyle\Upsilon_{u}(t) & \geq \frac{g_{0}}{4}t^{2}\|u\|^{2}\nonumber-C'_{4}\frac{\epsilon}{q}t^{2}\|u\|^{2} -C_{5}\frac{|t|^{r}}{q}\|u\|^{r} - C_{3}\frac{1}{4}\epsilon t^{2}\|u\|^{2}\\ \nonumber
   &- C_{1}\left( \int_{B}|tu|^{a'r}dx\right)^{ \frac{1 }{a'} }\left(\int_{B}\exp (\alpha t a|u|^{\gamma})dx\right)^{\frac{1}{a}}\nonumber \\ &\geq \left (\frac{g_{0}}{4}-\epsilon(\frac{1}{4} C_{3}+C'_{4}\frac{1}{q}) \right)\|tu\|^{2}-C_{4}\frac{|t|^{r}}{q}\|u\|^{r}- \left(\int_{B}\exp \big(\alpha  a\|tu\|^{\gamma}\big(\frac{|u|}{\|u\|}\big)^{\gamma}\big)dx\right)^{\frac{1}{a}} C_{5}\|tu\|^{r}\nonumber
   \end{align}
  By (\ref{eq:1.3}), the last integral is finte provided $t>0$ is chosen small enough such that  $\displaystyle\alpha a\|tu\|^{\gamma}\leq \alpha_{\beta}$. Then,
\begin{align*}
\displaystyle\Upsilon_{u}(t) & \geq  \left (\frac{g_{0}}{4}-\epsilon(\frac{1}{4} C_{3}+C'_{4}\frac{1}{q}) \right)\|tu\|^{2}- C_{6}\|tu\|^{r}~~\mbox{with}~~\alpha a\|tu\|^{\gamma}\leq \alpha_{\beta}~~\mbox{and}~~\alpha >\alpha_{0}\end{align*}
holds. Choosing  $\epsilon > 0$ such that $\displaystyle\frac{g_{0}}{4}-\epsilon(\frac{1}{4} C_{3}+C'_{4}\frac{1}{q}) > 0$ and since $ r > 4$, we obtain,\begin{equation}\label{eq:3.4}  \displaystyle\Upsilon_{u}(t) > 0 ~~ \mbox{for small}~~ t>0 .\end{equation}

Now, from $(A_{2})$, we can derive that there exist $C_{5}, C_{6} > 0 $ such that
\begin{equation}\label{eq:3.5}
F (x,t)\geq  C_{5}|t|^{\theta}-C_{6}.
\end{equation}
Then, by using $(\ref{eq:1.6})$ and $(\ref{eq:3.5})$, we get
\begin{align*}
 \displaystyle\Upsilon_{u}(t)
  &= \mathcal{J}(tu ) \nonumber  \\
  & \leq \frac{g(1)}{2} t^{2}\|u\|^{2}+\frac{g(1)}{4}t^{4}\|u\|^{4}-C'_{5}|t|^{\theta}|u|_{\theta}^{\theta}
  - C_{6}|B|
\end{align*}
Since $\theta>4$, we obtain \begin{align}\displaystyle\Upsilon_{u}(t)\rightarrow -\infty~~\mbox{as}~~t\rightarrow+\infty .\end{align} Hence, from (\ref{eq:3.4}) and (\ref{eq:3.5}),  there exists at least one $t_{u} > 0$ such that $\Upsilon'_{u}(t_{u})=0$, i.e. $t_{u} u\in \mathcal{N}$.  \\

Now we will show the uniqueness of  $t_{u}$. Let $s>0$ such that $su\in \mathcal{N}$ and suppose that $s\neq t_{u}$. Without loss of generality, we can assume that $s> t_{u}$. So we have $\displaystyle\langle\mathcal{J}'(t_{u}u),  t_{u}u \rangle= 0$ and
 $\displaystyle\langle\mathcal{J}'(su), su \rangle= 0$, then

\begin{equation}\label{eq:3.7}
 \displaystyle \frac{g(\|su\|^{2})}{s^{2}\|u\|^{2}} =\frac{1}{\|u\|^{4}}\bigg(\int_{B}|su|^{q-4}u^{4}dx+\int_{B}\frac{f(x,su)}{(su)^{3}}u^{4} dx\bigg)
\end{equation}

\begin{equation} \label{eq:3.8}
  \displaystyle \frac{g(\|t_{u}u\|^{2})}{t_{u}^{2}\|u\|^{2}} =\frac{1}{\|u\|^{4}}\bigg(\int_{B}|t_{u}u|^{q-4}u^{4}dx+\int_{B}\frac{f(x,t_{u}u)}{(t_{u}u)^{3}}u^{4} dx\bigg)\cdot
\end{equation}

Combining (\ref{eq:3.7}), (\ref{eq:3.8}), we get
 $$\displaystyle \frac{g(\|su\|^{2})}{s^{2}\|u\|^{2}}-\displaystyle \frac{g(\|t_{u}u\|^{2})}{t_{u}^{2}\|u\|^{2}} =\frac{1}{\|u\|^{4}}\bigg(\int_{B}\big((|su|^{q-4}-|t_{u}u|^{q-4})u^{4}+(\frac{f(x,su)}{(su)^{3}}-\frac{f(x,t_{u}u)}{(t_{u}u)^{3}})u^{4}\big)dx\bigg)\cdot$$
 Clearly, according to $(A_{3})$, Remark \ref{rem1} and $(G_{2})$, the left-hand side of the last equality is negative for
$t_u>s$ while the right-hand side is positive, which is a
	contradiction.
  This contradict the fact  that $s> t_{u}$. The case $
t_{u} > s>0$ is similar and we omit it. Then, $s= t_{u}$.\\
$(ii)$ Follows from $(i)$ , since $\mathcal{J}(t_{u}u)=\displaystyle\max_{t\geq 0}\Upsilon_{u}(t)$. \hfill $\Box$
\begin{lemma}\label{lem2.2}
	Assume that  $(A_1)-(A_4)$ hold. Then for any $u \in E$ with $u \neq 0$ such that
	$\langle \mathcal{J}'(u),u\rangle\leq0$, the unique maximum point  of  $\Upsilon_{u}$ on $\mathbb{R}_+ $ satisfies $0< t_u\leq 1$.
\end{lemma}

\noindent Proof:\\
	
	Since $t_u u\in\mathcal{N}$, we have
	\begin{equation} \label{eq:3.9}
	\begin{aligned}
\displaystyle \frac{g(\|t_{u}u\|^{2})}{t_{u}^{2}\|u\|^{2}} =\frac{1}{\|u\|^{4}}\bigg(\int_{B}|t_{u}u|^{q-4}u^{4}dx+\int_{B}\frac{f(x,t_{u}u)}{(t_{u}u)^{3}}u^{4} dx\bigg)\cdot
	\end{aligned}
	\end{equation}
	Furthermore, since $\langle \mathcal{J}'(u),u\rangle\leq0$, we have
	$$  \displaystyle \frac{g(\|u\|^{2})}{\|u\|^{2}} \leq\frac{1}{\|u\|^{4}}\bigg(\int_{B}|u|^{q-4}u^{4}dx+\int_{B}\frac{f(x,u)}{(u)^{3}}u^{4} dx\bigg).$$
	Then by  \eqref{eq:3.9}, we have
	\begin{equation} \label{eq:3.10}
	\begin{aligned}
\displaystyle \frac{g(\|t_{u}u\|^{2})}{t_{u}^{2}\|u\|^{2}}	- \displaystyle \frac{g(\|u\|^{2})}{\|u\|^{2}}
	\geq \frac{1}{\|u\|^{4}} \int_{B}\Big((|t_{u}u|^{q-4}-|u|^{q-4}+\frac{f(x,t_u u)}{(t_u u)^{3}}-\frac{f(x, u)}{u^{3}})\Big)u^4\,dx.
	\end{aligned}
	\end{equation}
	Obviously, from $(G_{2})$ the left hand side of \eqref{eq:3.10} is negative for
	$t_u>1$ whereas the right hand side is positive, which is a
	contradiction. Therefore $0<t_u\leq 1$.\hfill $\Box$\\
In the sequel, we prove that sequences in $\mathcal{N}$ cannot converge to $0$.
\begin{lemma}\label{lem3}
 For all $u\in \mathcal{N}$,
 \begin{itemize}
   \item [$(i)$] there exists $\kappa>0$ such that\\
   $ \|u\|  \geq \kappa ;$
   \item[$(ii)$] $\mathcal{J}(u) \geq \displaystyle(\frac{1}{4}-\frac{1}{q})g_{0}\|u\|^{2}$
 \end{itemize}
 \end{lemma}
  Proof. $(i)$
We argue by contradiction.  Suppose that there exists a sequence $\{u_{n}\} \subset \mathcal{N} $
such that $u_{n}\rightarrow 0$ in $ \mathbf{W}.$
Since $\{u_{n}\} \subset \mathcal{N}$, then  $\displaystyle\langle\mathcal{J}'(u_{n}) ,u_{n} \rangle =0$.
Hence, it follows from (\ref{eq:3.2}), (\ref{eq:3.3}) and the radial Lemma \ref{lemr} that

\begin{align}\label{eq:3.11}
g_{0}\|u_{n}\|^{2}<g(\|u_{n}\|^{2})\|u_{n}\|^{2}&=\int_{B} |u|^{q}dx+\int_{B}f(x,u_{n})u_{n}dx \\
 & \leq 2\epsilon  \int_{B}|u_{n}|^{2}dx + C_{1}\int_{B}|u_{n}|^{r}dx +C'_{1}\int_{B}|u_{n}|^{r}\exp(\alpha |u_{n}|^{\gamma} )dx\nonumber\\
 & \leq \epsilon  C_{6} \|u_{n}\|^{2} +C_{7}\|u\|^{r} + C_{1}\int_{B}|u_{n}|^{r}\exp(\alpha |u_{n}|^{\gamma} )dx\nonumber
\end{align}
Let $a>1$ with $\frac{1}{a}+\frac{1}{a'}=1$. Since $u_{n}\rightarrow 0\mbox{ in }~~\mathbf{W}$,
for $n$ large enough, we get
$\displaystyle\|u_{n}\|\leq(\frac{\alpha_{ \beta}}{\alpha a})^{\frac{1}{\gamma}}$. From H\"{o}lder inequality, (\ref{eq:1.5}) and again the radial Lemma \ref{lemr},
 we have
\begin{align*}
  \int_{B}|u_{n}|^{r}\exp(\alpha |u_{n}|^{\gamma} )dx&\leq
  \left( \int_{B}|u_{n}|^{ra'}dx\right)^{\frac{1}{a'}}
  \left(\int_{B}\exp \big(\alpha  a\|u\|^{\gamma}\big(\frac{|u|}{\|u\|}\big)^{\gamma}\big)dx\right)^{\frac{1}{a}} \\
  &\leq C_{7} \left( \int_{B}|u_{n}|^{ra'}dx\right)^{\frac{1}{a'}}\leq C_{8} \|u_{n}\|^{r}
\end{align*}
Combining (\ref{eq:3.11}) with the last inequality, for $n$ large enough, we obtain

\begin{equation}\label{eq:3.12}
g_{0} \|u_{n}\|^{2}\leq \epsilon C_{6} \|u_{n}\|^{2}
 + C_{8} \|u_{n}\|^{r}
\end{equation}
Choose suitable $\epsilon > 0 $ such that  $g_{0}-\epsilon C_{6} > 0$.
Since  $2 < r$, then (\ref{eq:3.12}) contradicts the fact that $u_{n}\rightarrow 0\mbox{ in }\mathbf{W}$.\\
$(ii)$ Given $u \in \mathcal{N}$, by the deﬁnition of $\mathcal{N}$, $(\ref{eq:1.8})$ and $(A_{3})$, we obtain
\begin{align*}
  \mathcal{J}(u) &= \mathcal{J}(u)-\frac{1}{q}\langle\mathcal{J}'(u), u\rangle \\
  &= \frac{1}{2}G(\|u\|^{2})-\frac{1}{q}g(\|u\|^{2})\|u\|^{2}
  + \big(\int_{B}\frac{1}{q}f(x,u)u-F(x,u)dx\big)\\
  &\geq (\frac{1}{4}-\frac{1}{q})g_{0}\|u\|^{2}
\end{align*}
Lemma \ref{lem3} implies that $\mathcal{J}(u) > 0$  for all $u\in \mathcal{N}$.\hfill $\Box$\\
As a consequence, $\mathcal{J}$ is bounded by below in $\mathcal{N}$, and therefore
$\displaystyle m:=\inf_{u\in \mathcal{N}} \mathcal{J}(u) $
 is well-deﬁned.\\ In the following lemma we prove that if the minimum of $\mathcal{J}$ on $\mathcal{N}$  is realized at some  $u\in\mathcal{N}$, then $u$  is a critical point of $\mathcal{J}$.

\begin{lemma}\label{lem5}
If $ u_{0}\in \mathcal{N} $ satisfies $\mathcal{J}(u_{0})=m$, then $\displaystyle\mathcal{J}'(u_{0})=0.$
\end{lemma}
\noindent Proof:We argue by contradiction. We assume that $\displaystyle\mathcal{J}'(u_{0})\neq 0$. By the continuity of $\mathcal{J}'_{\lambda}$,
there exist $\iota, \delta\geq 0$ such that
\begin{equation}\label{eq:3.13}
\displaystyle \|\mathcal{J}'_{\lambda}(v)\|_{\mathbf{W}^{\ast}}\geq\iota \mbox{ for all }~~v ~~\mbox{such that}~~ \|v-u_{0}\|\leq \delta.
\end{equation}
 Let $\displaystyle D=\left(1-\tau,1+\tau\right)\subset \mathbb{R}$ with $\displaystyle \tau \in(0, \frac{\delta}{4\|u_{0}\|})$ and define $h:D\rightarrow\mathbf{W}$, by

$$\displaystyle h(\rho)=\rho u_{0},  \rho \in D$$
By virtue of $u_{0} \in \mathcal{N}$, $\mathcal{J}(u_{0})=m$ and Lemma \ref{lem1}, it is easy to see that
\begin{equation}\label{eq:3.14}
\displaystyle\bar{m}:=\max_{\partial D} \mathcal{J}\circ h<m~~ \mbox{and}~~ \mathcal{J}(h(\rho))<m,~~\forall ~~\rho\neq 1.
\end{equation}
 Let $\epsilon:=\min\{\frac{m-\bar{m}}{2}, \frac{\iota\delta}{16}\}$, $S_{r}:=B(u_{0},r),r\geq0$
and $\displaystyle\mathcal{J}^{a}:=\mathcal{J}^{-1}(]-\infty,a]).$
 According to the quantitative deformation Lemma $[\cite{Wi}, \mbox{ Lemma }2.3]$,
there exists a deformation $\eta \in C\left(\mathbf{W}, \mathbf{W}\right)$ such that:
 \begin{itemize}
   \item [$(1)$] $\eta( v)=v,$ if $v\not\in \mathcal{J}^{-1}([ m-\epsilon,m+\epsilon])\cap S_{\delta}$
   \item[$(2)$] $\eta\left( \mathcal{J}^{m+\epsilon }\cap S_{\frac{\delta}{2}}\right)\subset \mathcal{J}^{m-\epsilon}$,
    \item[$(3)$] $\mathcal{J}(\eta( v))\leq \mathcal{J}(v)$, for all $v\in  \mathbf{W}. $
 \end{itemize}
By lemma \ref{lem1} $(ii)$, we have $\mathcal{J}( h(\rho))\leq m$. In addition, we have,$$\|h(\rho)-u_{0}\|=\|(\rho-1)u_{0}\|\leq \frac{\delta}{4},~~\forall \rho\in D\cdot$$ Then $h(\rho)\in S_{\frac{\delta}{2}}$ for $\rho\in \bar{D}$.
Therefore, it follows from $(2)$ that
\begin{equation}\label{eq:3.15}
\max_{\rho\in \bar{D}}\mathcal{J}(\eta( h(\rho))\leq m-\epsilon.
\end{equation}
In the sequel, we will prove that $\eta( h(D))\cap\mathcal{N}$ is nonempty. In such case, due to the definition of $m$, this contradicts (\ref{eq:3.15}).
 To do this, we first define
$$\bar{h}(\rho):=\eta( h(\rho)),$$
\begin{align*}
  \Upsilon_{0}(\rho)& = \langle \mathcal{J}'(h(\rho)), u_{0} \rangle,
\end{align*}
and
$$\Upsilon_{1}(\rho):=(\frac{1}{\rho}\langle \mathcal{J}'(\bar{h}(\rho),(\bar{h}(\rho)) \rangle.$$
We have that for $\rho\in \overline{D}$, $\mathcal{J}(h(\rho))\leq \overline{m}<m-\varepsilon$. Then, $\bar{h}(\rho)=\eta (h(\rho))=\rho u_{0}$. Hence,
\begin{align}\label{eq:3.16}
  \Upsilon_{0}(\rho)& =  \Upsilon_{1}(\rho), \forall \rho \in \overline{D}
\end{align}
On one hand, we have that $\rho=1$ is the unique critical point of $\Upsilon_{0}$. So by degree theory, we get that $d^{0}(\Upsilon_{0},\mathcal{J},0)=1$. On the other hand, from (\ref{eq:3.16}), we deduce that  $d^{0}(\Upsilon_{1},\mathcal{J},0)=1$. Consequently, there exists $\overline{\rho}\in D$ such that $\overline{h}(\overline{\rho})\in \mathcal{N}$. This implies that
$$m\leq \mathcal{J}(\overline{h}(\overline{\rho}))=\mathcal{J}(\eta(h(\overline{\rho})).$$
This contradicts (\ref{eq:3.15}) and finish the proof of the Lemma.\hfill $\Box$\\
\section{Auxiliary problem}
In this section, in order to prove our existence  result , we consider
the auxiliary problem \begin{equation}\label{pau}
 \displaystyle (P_{a})~~~~\left\{
\begin{array}{rclll}
g\big(\int_{B}w(x) |\Delta u|^{2}\big) \Delta (w (x) \Delta u) &=&  \displaystyle \vert u\vert^{p-2}u& \mbox{in} & B \\\\
\displaystyle u&=&\displaystyle\frac{\partial u}{\partial n}=0 &\mbox{on }&  \partial B,
\end{array}
\right.
\end{equation}
where $p$ is the constant that appear in the hypothesis $(A_5)$. The energy $J_p$  associated to problem \eqref{pau} is given by
$$J_p(u) :=\frac{1}{2}G(\|u\|^{2}) -\frac{1}{p}\int_{B} \vert u\vert^p dx . $$ We introduce  the Nehari manifold associeted to $J_p$ that is
$$\mathcal{N}_p:= \{ u \in \mathbf{W}, u \neq 0\; \mbox{and}\; \langle J_p'(u),u \rangle = 0\}.$$
Let $m_p = \displaystyle\inf_{\mathcal{N}_p} J_p(u)>0$, we have the following results for $J_{p}$.
\begin{lemma}\label{ljp1} Given $u\in \mathbf{W}, u\neq 0$, there exists a unique $t>0$ such that $t u\in\mathcal{N}_{p}$. In addition, $t$ satisfies \begin{equation}\label {JP1}J_{p}(tu)=\displaystyle \max_{s\geq0}J_{p}(su).\end{equation}
\end{lemma}
\noindent Proof: Let $$\displaystyle\gamma(s) = J_{p}(su )=\frac{1}{2}G(s^{2}\|u\|^{2})-\frac{s^{p}}{p}|u|^{p}_{p}\geq \frac{g_{0}}{4}s^{2}\|u\|^{2}-\frac{s^{p}}{p}|u|^{p}_{p}$$ for $s>0$. Since $p>4$, we have that $\gamma(s)>0$ for $s>0$ small enough and $\gamma(s)\rightarrow -\infty~~\mbox{as}~~~s\rightarrow -\infty$. Hence, there exists $t>0$ satisfying (\ref{JP1}). In particular, $tu\in \mathcal{N}_{p}$. We can proceed as in Lemma \ref{lem1} to prove that  $t$ is unique .\hfill $\Box$\\

As a consequence, we have
\begin{corollary} Let $u\in \mathbf{W}, u\neq 0$. Then  $u\in \mathcal{N}_{p}$ if and only if $J_{p}(tu)=\displaystyle \max_{s\geq0}J_{p}(su)$.
\end{corollary}
Also, it is easy to proof the following Lemmas:
\begin{lemma}\label{ljp2}
 For all $u\in \mathcal{N}_{p}$,
 \begin{itemize}
   \item [$(i)$] there exists $\kappa_{0}>0$ such that\\
   $ \|u\|  \geq \kappa_{0} ;$
   \item[$(ii)$] $\mathcal{J}_{p}(u) \geq (\frac{1}{4}-\frac{1}{p})|u|^{p}_{p}$
 \end{itemize}
 \end{lemma}
\begin{lemma} \label{ljp3}There exists $w_{p} \in \mathcal{N}_{p}$ such that $J_{p}(w_{p}) = m_{p}$.
\end{lemma}
\noindent Proof:
	Let sequence $(w_n) \subset \mathcal{N}_{p} $ satisfy $\displaystyle \lim_{n \rightarrow +\infty} J_{p}(w_n) = m_{p}$. It is clearly that $(w_n)$ is bounded by Lemma \ref{ljp2}. Then, up to a subsequence, there exists $w_{p} \in \mathbf{W}$ such that
	
	\begin{equation}\label{eq:4.3}\begin{array}{ll}
	w_{n} \rightharpoonup w_{p}~~~~&\mbox{in}~~\mathbf{W},\\
	w_{n} \rightarrow w_{p}~&\mbox{in}~~L^{q}(B),~~\forall q\geq 2,\\
	w_{n} \rightarrow w_{p} ~~&\mbox{a.e. in }~~B.
	\end{array}
	\end{equation}
	We claim that $w_{p}\ne 0$. Suppose, by contradiction, $w_{p}= 0$. From the definition
	of $\mathcal{N}_{p}$ and \eqref{eq:4.3} , we have that $ \displaystyle\lim_{n \rightarrow +\infty} \Vert w_n \Vert^2 =0$, which contradicts Lemma \ref{ljp2}. Hence, $w_{p}\ne 0$ .
	
	From the continuity of $g$, the lower semi continuity of norm and \eqref{eq:4.3}, it follows that
 	\begin{equation}\label{eq:4.4}
 g(\Vert w_{p}\Vert^2)\Vert w_{p}\Vert^2 \leq  \liminf_{n \rightarrow +\infty}g(\Vert w_n\Vert^2 )\Vert w_n\Vert^2
 	\end{equation}	

 On the other hand, by using $\langle J'(w_n),w_n\rangle =0$ and (\ref{eq:4.3}), we have

 	\begin{equation}\label{eq:4.5} \liminf_{n \rightarrow +\infty}g(\Vert w_n\Vert^2 )\Vert w_n\Vert^2 = \liminf_{n \rightarrow +\infty} \int_{B}| w_n|^{p} dx=\int_{B} |w_{p}|^{p} dx.\end{equation}
 From \eqref{eq:4.4} and \eqref{eq:4.5} we deduce that $\langle J_{p}'(w_{p}),w_{p}\rangle \leq 0$.
 Then, as in Lemma \ref{lem2.2} this implies that there exists $s_u \in (0, 1] $ such that $s_uw_{p}\in\mathcal{N}_{p}$. Thus, by the lower semi continuity of norm, (\ref{eq:1.8}) and (\ref{eq:4.3}), we get that

$$\begin{array}{rclll}
\displaystyle m_{p} \leq J_{p}(s_uw_{p})&=&J(s_uw_{p}) -\displaystyle\frac{1}{4} \langle J_{p}'(s_uw_{p}),s_uw_{p}\rangle\\&=&\displaystyle\frac{1}{2}G(\|s_{u}w_{p}\|^{2})-\frac{1}{4}g(\|s_{u}w_{p}\|^{2})\|s_{u}w_{p}\|^{2}+\displaystyle\big(\frac{1}{4}-\frac{1}{p}\big)s^{p}_{u}\int_{B}|w_{p}|^{p}dx
\\ &\leq& J_{p}(w_{p}) -\displaystyle\frac{1}{4} \langle J_{p}'(w_{p}),w_{p}\rangle\\
&=&\displaystyle\frac{1}{2}G(\|w_{p}\|^{2}) -\frac{1}{p}\int_{B}|w_{p}|^{p}dx-\frac{1}{4}g(\|w_{p}\|^{2})\|w_{p}\|^{2}+\frac{1}{4}\int_{B}|w_{p}|^{p}dx\\
&\leq&\displaystyle\liminf_{n\rightarrow +\infty}\Big[\displaystyle\frac{1}{2}G(\|w_{n}\|^{2}) -\frac{1}{p}\int_{B}|w_{n}|^{p}dx\bigg]\\&-&\displaystyle\liminf_{n\rightarrow +\infty}\Big[\frac{1}{4}g(\|w_{n}\|^{2})\|w_{p}\|^{2}+\frac{1}{4}\int_{B}|w_{n}|^{p}dx\Big]\\
&\leq&\displaystyle \liminf_{n\rightarrow+\infty}\big[J_{p}(w_n) -\displaystyle\frac{1}{4} \langle J_{p}'(w_n), w_n\rangle\big] = m_{p}.
\end{array}$$	
Therefore, we get that $J_{p}(w_{p})=m_{p}$, which is the desired conclusion.\hfill $\Box$\\
\section{Proof of Theorem 1.2}

Now we will get an essential estimate for level $m$. This will be a useful tool
to obtain an adequate bound on the norm of a minimizing sequence for $m$ in $\mathcal{N}$.
	\begin{lemma}\label{lem9}
	Assume that $(A_1)-(A_5)$ and \eqref{eq:1.9} are satisfied. It holds that
		\begin{equation} \label{eq:5.1} m \leq g_{0}\frac{q-4}{4q}\bigg(\frac{ \alpha_{\beta}}{2 (\alpha_0 + \delta)}\bigg)^{1-\beta}.\end{equation}
\end{lemma}

\noindent Proof:
From Lemma \ref{ljp3}, there exists $w_{p} \in \mathcal{N}_p$ such that $J_p(w_{p}) = m_p$ and $J_p'(w_{p})=0.$ Consequently, using (\ref{eq:1.8}) we get	
\begin{equation}\label{eq:5.2}
\frac{1}{2}G(\Vert w_{p}\Vert^2) - \frac{1}{p}\int_{B} \vert w_{p}\vert^p ~dx =  m_p
\end{equation}
and
	\begin{equation} \label{eq:5.3}
g(\Vert w_{p}\Vert^2)\Vert w_{p}\Vert^2 = \int_{B} \vert w_{p}\vert^p ~dx.
\end{equation}
Note that by using \eqref{eq:5.2}, \eqref{eq:5.3} and the fact that $p>q>4$, we have
$$(\frac{1}{q}- \frac{1}{p})|w_{p}|_p^p  =\frac{1}{q}g(\|w_{p}\|^{2})\|w_{p}\|^{2}- \frac{1}{2}G(\|w_{p}\|^{2})+m_{p}\leq m_p.$$
So,
\begin{equation}\label{eq:5.4}
 |w_{p}|^{p}_{p}<\frac{pq}{p-q}m_{p}
\end{equation}
According to $(A_5)$ and \eqref{eq:5.3}, we have $\langle \mathcal{J}'(w_{p}),w_{p}\rangle \leq 0$ which, with lemma \ref{lem2.2}, gives that there exists a unique $s \in (0, 1) $ such that $sw_{p}\in\mathcal{N}$. Using $(A_5)$, \eqref{eq:5.2}, \eqref{eq:5.3}, \eqref{eq:1.7} and \eqref{eq:5.4}, we obtain
$$\begin{array}{rclll}
\displaystyle m &\leq& \mathcal{J}(sw_{p})
\\  &\leq& \displaystyle\frac{g(1)s^2}{2}\|w_{p}\|^2+\frac{g(1)s^4}{4}\|w_{p}\|^4
\displaystyle - \frac{C_p s^p}{p}\vert w_{p}\vert_p^p \\  &\leq& \displaystyle\frac{g(1)s^2}{2}\|w_{p}\|^2+\frac{g(1)s^2}{4}\|w_{p}\|^4
\displaystyle - \frac{C_p s^p}{p}\vert w_{p}\vert_p^p \\ &\leq& \displaystyle\frac{g(1)s^2}{2g_{0}}|w_{p}|^{p}_{p}+\frac{g(1)s^2}{4g^{2}_{0}}|w_{p}|^{2p}_{p}
\displaystyle - \frac{C_p s^p}{p}\vert w_{p}\vert_p^p \\&=& \displaystyle\bigg((\frac{g(1)}{2g_{0}}+\frac{g(1)}{4g^{2}_{0}}|w_{p}|^{p}_{p})s^{2} -  \frac{C_p s^p}{p}\bigg)\vert w_{p}\vert_p^p  \\
 &\leq&\displaystyle \max_{\xi >0}\bigg(\big(\frac{g(1)}{2g_{0}}+\frac{g(1)}{4g^{2}_{0}}|w_{p}|^{p}_{p}\big)\xi^{2} -  \frac{C_p \xi^p}{p}\bigg)\vert w_{p}\vert_p^p\\ &\leq&\displaystyle \max_{\xi >0}\bigg(\big(\frac{g(1)}{2g_{0}}+\frac{g(1)}{4g^{2}_{0}}\frac{pq}{p-q}m_{p}\big)\xi^{2} -  \frac{C_p \xi^p}{p}\bigg)\vert w_{p}\vert_p^p
\end{array}$$
By some simple algebraic calculations, we get
\begin{equation}\label{eq:5.5}
m \leq\tau (\frac{2\tau}{C_{p}})^{\frac{2}{p-2}}\big(\frac{p-2}{p}\big )|w_{p}|_p^p.
\end{equation}
Thus, by using \eqref{eq:5.5}, we obtain
\begin{equation}\label{eq:5.6}
m < \tau(\frac{2\tau}{C_{p}})^{\frac{2}{p-2}}\big(\frac{q(p-2)}{p-q}\big )m_p.
\end{equation}
Therefore, by \eqref{eq:1.9} and \eqref{eq:5.6}, we get that \eqref{eq:5.1} is valid.\hfill $\Box$\\

The result below gives us some compactness properties of minimising sequences.
	\begin{lemma}\label{lemma9}
	If $(u_n) \subset \mathcal{N} $ is a minimizing sequence for $m$, then there exists $u \in \mathbf{W}$ such that
	
	$$\int_{B}f(x,u_{n})u_n dx \rightarrow \int_{B}f(x,u)u dx $$
	and
	$$\int_{B}F(x,u_{n}) dx \rightarrow \int_{B}F(x,u)dx. $$
\end{lemma}

\noindent Proof:
	We must prove the first limit, since the second one is analogous. For this, we use (\ref{imp}) and introduce the following function  $k(u_n(x))$ given by $$k(u_n(x)):=\e
	| u_n|^2\,dx + C|u_n|^{q} \exp(\alpha |u_n|^{\gamma}).$$   It's clear  that  is sufficient to prove that $k(u_n(x))$ is convergent in $L^1(B)$.  We have \begin{equation*}\label{eq4.1}
	\int_{B}f(x, u_n)~ u_n dx \leq \e
	\int_{B}| u_n|^2\,dx + C \int_{B}|u_n|^{q} \exp(\alpha |u_n|^{\gamma})dx =\int_{B} k(u_n(x))~dx, \quad \mbox{for all}\; \alpha > \alpha_0\, \mbox{and}\, q>2.
	\end{equation*}
First note that from Sobolev embedding, we have
	\begin{equation}\label{eq:5.7}
	|u_{n}|^q \rightarrow |u|^q~~ \mbox{in}~~ L^{2}(B).
	\end{equation}	
	On the other hand, by $(A_{2})$ and (\ref{eq:1.7'}), we obtain that
		\begin{equation} \label{eq:5.8}
	\begin{aligned}
	\displaystyle m=\limsup_{n \rightarrow +\infty}  \mathcal{J}(u_n)
	& = \limsup_{n \rightarrow +\infty} \big(\mathcal{J}(u_n)-\frac{1}{q}\langle \mathcal{J}'(u_n),u_n\rangle \big) \\
	& =\limsup_{n \rightarrow +\infty} \big(  \frac{1}{2}G(\| u_n\|^{2})-\frac{1}{q} g(\| u_n\|^{2})\| u_n\|^{2} +\frac{1}{q}\int_{B}\big(f(x, u_n) u_n- q F(x, u_n)\big)dx\big) \\
	& >  (\frac{1}{4}-\frac{1}{q})\limsup_{n \rightarrow +\infty} g(\| u_n\|^{2}\|)\| u_n\|^{2}\\
    & >\frac{q-4}{4q}g_{0}\limsup_{n \rightarrow +\infty}\| u_n\|^{2}
	\end{aligned}
	\end{equation}
	which, together with Lemma \ref{lem9} leads to the following estimation  $\displaystyle \limsup_{n \rightarrow +\infty} \|u_n\|^\gamma < \frac{\alpha_{\beta}}{2(\alpha_0+\delta)}\cdot$
	
    Now choosing $\alpha= \alpha_0+\delta,~~\delta>0$, we have  that
	\begin{equation}\label{eq:5.9}
	\int_{B} \exp(2\alpha |u_n|^{\gamma})dx \leq \int_{B}\exp\big(2(\alpha_{0} + \delta )\|u_{n}\|^{\gamma}(\frac{|u_{n}|}{\|u_{n}\|})^{\gamma}\big)dx\leq \int_{B}\exp\big(\alpha_{\beta}(\frac{|u_{n}|}{\|u_{n}\|})^{\gamma}\big)dx.
	\end{equation}
	Then it follows by (\ref{eq:1.3}) that there is $M > 0$ such that
		\begin{equation}\label{eq:5.10}
 \int_{B} \exp(2 \alpha |u_n|^{\gamma})dx \leq M.
	\end{equation}
	Since
	\begin{equation}\label{eq:5.11}
	\exp(\alpha |u_n|^{\gamma}) \rightarrow \exp(\alpha |u|^{\gamma})\:  \mbox{a.e in }~~B,
	\end{equation}
	from \eqref{eq:5.9} and [\cite{Wil}, Lemma 4.8], we get that
	\begin{equation}\label{eq:5.12}
	\exp(\alpha |u_n|^{\gamma}) \rightharpoonup \exp(\alpha |u|^{\gamma})\:  \mbox{ in }~~L^{2}(B).
	\end{equation}
	Now using \eqref{eq:5.7}, \eqref{eq:5.12} and [\cite{Wil}, Lemma 4.8] again, we conclude that
	\begin{equation}\label{eq:5.13}
\int_{B} f(x, u_n) ~u_n dx \rightarrow \int_{B} f(x, u)~ u dx.
	\end{equation}\hfill $\Box$\\

In the following, we give an additional important result that will be used to prove our main result.
\begin{lemma}\label{lem11} Assume that the conditions ($A_{1}$), ($A_{2}$) and ($A_{3}$) are satisfied. Then, for each $x\in B$, we have
	$$tf(x, t) -q F(x, t)~~ \mbox{is increasing for} ~ t> 0 ~ \mbox{ and decreasing for} ~ t< 0.  $$
	In particular, $tf(x, t) - q F(x, t)> 0~~\mbox{for all} ~~(x,t) \in B \times \R\setminus \{0\}.$
\end{lemma}
\noindent Proof: Assume that $0<t<s$. For each $x\in B$, we have
$$\begin{array}{rlll}
\displaystyle
tf(x,t)-qF(x,t)&=& \displaystyle\frac{f(x,t)}{t^{q-1}}t^{q}-qF(x,s)+q\int^{s}_{t}f(x,\nu)d\nu \\
&<&\displaystyle \frac{f(x,t)}{s^{q-1}}t^{q}-qF(x,s)+\frac{f(x,s)}{s^{q-1}}(s^{q}-t^{q})\\

&=&sf(x,s)-q F(x,s)\cdot
\end{array}$$
The proof in the case $t<s<0$ is similar.\\
The assertion $tf(x, t) - q F(x, t)> 0~~\mbox{for all} ~~(x,t) \in B \times \R\setminus \{0\}$ comes from ($A_{2}$).\hfill $\Box$\\

By the following lemma, we prove that the minimum of $\mathcal{J}$ on $\mathcal{N}$ is achieved in some $w_{0}\in \mathcal{N}$.

\begin{lemma} \label{lem12}There exists $w_{0} \in \mathcal{N}$ such that $\mathcal{J}(w_{0}) = m$ .
\end{lemma}
\noindent Proof:
	Let sequence $(w_n) \subset \mathcal{N} $ satisfying  $\displaystyle \lim_{n \rightarrow +\infty}\mathcal{J} (w_n) = m$. It is clearly that $(w_n)$ is bounded by Lemma \ref{ljp2}. Then, up to a subsequence, there exists $w_{0} \in E$ such that
	
	\begin{equation}\label{eq:5.14}\begin{array}{ll}
	w_{n} \rightharpoonup w_{0}~~~~&\mbox{in}~~E,\\
	w_{n} \rightarrow w_{0}~&\mbox{in}~~L^{q}(B),~~\forall q\geq 2,\\
	w_{n} \rightarrow w_{0} ~~&\mbox{a.e. in }~~B.
	\end{array}
	\end{equation}
	We claim that $w_{0}\ne 0$. Suppose, by contradiction, $w_{0}= 0$. From the definition
	of $\mathcal{N}$ and \eqref{eq:5.14} , we have that $ \displaystyle\lim_{n \rightarrow +\infty} \Vert w_n \Vert^2 =0$, which contradicts Lemma \ref{lem3}. Hence, $w_{0}\ne 0$ .
	
	From the lower semi continuity of norm, the continuity of $g$ and \eqref{eq:5.14}, it follows that
 	\begin{equation}\label{eq:5.15}
 	g(\Vert w_{0}\Vert^2)\Vert w_{0}\Vert^2 - \lim_{n\rightarrow \infty}\int_{B}|w_{0}|^{q}~~dx\leq  \liminf_{n \rightarrow +\infty}\big(g(\Vert w_n\Vert^2 )\Vert w_n\Vert^2-\int_{B}|w_{n}|^{q}~~dx\big)\cdot
 	\end{equation}	

 On the other hand, by using $\langle \mathcal{J}'(w_n),w_n\rangle =0$ and (\ref{eq:5.14}), we have

 	\begin{equation}\label{eq:5.16} \liminf_{n \rightarrow +\infty}g(\Vert w_n\Vert^2 )\Vert w_n\Vert^2 =\liminf_{n \rightarrow +\infty} \int_{B}(f(x, w_n) w_n+|w_{n}|^{q}) dx=\int_{B}(f(x, w_{0}) w_{0}+|w_{0}|^{q}) dx.\end{equation}
 From \eqref{eq:5.15} and \eqref{eq:5.16} we deduce that $\langle \mathcal{J}'(w_{0}),w_{0}\rangle \leq 0$.
 Then, as in Lemma \ref{lem2.2} this implies that there exists $s\in (0, 1] $ such that $sw_{0}\in\mathcal{N}$. Thus, by the lower semi continuity of norm, \eqref{eq:1.7}, Lemma \ref{lem11} and Lemma \ref{lemma9}, we get that

\begin{align}\nonumber\begin{array}{rclll}
\displaystyle m \leq \mathcal{J}(sw_{0})&=&\mathcal{J}(sw_{0}) -\displaystyle\frac{1}{q} \langle \mathcal{J}'(sw_{0}),sw_{0}\rangle\\&=&\displaystyle\frac{1}{2}G(\|sw_{0}\|^{2}) -\frac{1}{q}g(\|sw_{0}\|^{2})\|sw_{0}\|^{2}\\&+&\displaystyle\frac{1}{q}\int_{B}\big(f(x,sw_{0})sw_{0}-qF(x,sw_{0})\big)dx
\\ &<& \displaystyle\frac{1}{2}G(\|w_{0}\|^{2}) -\frac{1}{q}g(\|w_{0}\|^{2})\|w_{0}\|^{2}\\&+&\displaystyle\frac{1}{q}\int_{B}\big(f(x,w_{0})w_{0}-qF(x,w_{0})\big)dx \\
&\leq&\displaystyle\liminf_{n\rightarrow +\infty}\Big[\displaystyle\frac{1}{2}G(\|w_{n}\|^{2}) -\int_{B}F(x,w_{n})dx\big]\\&-&\displaystyle\liminf_{n\rightarrow +\infty}\Big[\frac{1}{q}g(\|w_{n}\|^{2})\|w_{n}\|^{2}-\frac{1}{q}\int_{B}f(x,w_{n})w_{n}dx\Big]\\
&\leq&\displaystyle \liminf_{n\rightarrow+\infty}\big[\mathcal{J}(w_n) -\displaystyle\frac{1}{q} \langle \mathcal{J}'(w_n), w_n\rangle\big] = m.
\end{array}\end{align}	
Therefore, we get that $\mathcal{J}(sw_{0})=m$, which is the desired conclusion.\hfill $\Box$\\

\textbf{Proof of Theorem \ref{th1.3}. }From Lemma \ref{lem12} there exists $w_{0}$ such that $\mathcal{J}(w_{0})=m$. Now, by Lemma \ref{lem5}, we deduce that $\mathcal{J}'(w_{0})=0$.
So, $w_{0}$ is a solution to problem $(P)$ .\hfill $\Box$\\
\section*{ Statements and Declarations:}

 We declare that this manuscript is original, has not been published before and is not currently being considered for
publication elsewhere.

We confirm that the manuscript has been read and approved and that there are no other persons who satisfied the criteria for
authorship but are not listed.
 \section*{ Competing Interests:}

The authors declare that they have no known competing financial interests or personal relationships that could have appeared to influence the work reported in this paper.

\end{document}